\documentclass{elsart}
\usepackage{amsmath}
\usepackage{amssymb}
\usepackage{latexsym}
\usepackage[cmtip,arrow]{xy}
\usepackage{pb-diagram,pb-xy}

\journal{Journal of Pure and Applied Algebra}

\newcommand\badO[1] {\cB^0_{#1}}
\newcommand\bad[1] {\cB_{#1}}
\newcommand\badOV[1]{\cB^0_{#1}(X, V)}
\newcommand\badV[1]{\cB_{#1}(X, V)}
\newcommand\badL[1]{\cB_{#1}(X, L)}
\newcommand{\brref}[1]{(\ref{#1})}
\newcommand{\bslocus}[1]{{\rm Bs} |#1|}

\newcommand{\calo}{{\mathcal O}}
\newcommand{\cB}{\mathcal B}

\newcommand{\xiPT}{\tilde{\xiP}}
\newcommand{\xiP}{\xi'}
\newcommand{\cI}{\mathcal I}
\newcommand{\cL}{\mathcal L}

\newcommand{\red}[1]{X^{[#1]}_{(1,\dots, 1)}}
\newcommand{\tX}{\tilde{X}}
\newcommand{\ttau}{\tilde{\tau}}
\newcommand{\tx}{\tilde{x}}
\newcommand{\ty}{\tilde{y}}
\newcommand{\tz}{\tilde{z}}
\newcommand{\X}[1]{X^{[{#1}]}}

\newcommand{\comp}{\mathbb C}
\newcommand{\Pin}[1]{{\mathbb P}^{#1}}
\newcommand{\restrict}[2]{{#1}_{\mid _{#2}}}
\newcommand{\Proj}[1]{\mathbb{P}(#1)}
\newcommand{\oof}[2]{{\mathcal O}_{#1}({#2})}
\newcommand{\oofp}[2]{{\mathcal O}_{\mathbb{ P}^{#1}}({#2})}
\newcommand{\oofo}[1]{{\mathcal O}_{#1}}

\newcommand{\iofo}[1]{{\mathcal I}_{#1}}

\newcommand{\scroll}[1]{(\Proj{#1},\taut{#1})}
\newcommand {\xel} {(X, L)}
\newcommand\supp[1]{{\rm Supp}(#1)}
\newcommand{\taut}[1]{{\mathcal O}_{\mathbb{P}(#1)}(1)}

\newcommand{\Picof}[1]{\text{\rm Pic}(#1)}
\newcommand{\phiL}{\varphi_{L}}
\newcommand{\phiV}{\varphi_{V}}
\newcommand{\tensor}{\otimes}
\newcommand{\length}[1]{{\rm length}(#1)}

\newcommand{\ddui}[2]{\frac{\partial #1}{\partial u_{#2}}}
\begin{document}
\begin{frontmatter}
\title{Higher Order Bad Loci}
\author{Gian Mario Besana \corauthref{cor}}
\address{School of C.T.I - DePaul University\\
243 S. Wabash - Chicago IL 60604 USA}
\corauth[cor]{Corresponding author.}
\ead{gbesana@cti.depaul.edu}

\author{Sandra Di Rocco}
\address{K.T.H.Matematik \\ S-10044 Stockholm, Sweden}
\ead{dirocco@math.kth.se}

\author{Antonio Lanteri}
\address{Dip. di Matematica ``F. Enriques" -
Universit\`a degli Studi di Milano \\
 Via Saldini 50 - 20133 Milano, Italy}
\ead{lanteri@mat.unimi.it}

\begin{abstract} Zero-schemes on smooth complex projective varieties,
forcing all elements of ample and free linear systems to be
reducible are studied. Relationships among the minimal length of
such zero-schemes, the positivity of the line bundle associated
with the linear system, and the dimension of the variety are
established. Bad linear spaces are also investigated.
\end{abstract}
\end{frontmatter}
\section{Introduction}
\label{intro} Given a linear system on a smooth, complex,
projective variety $X$ with $\dim{X} \ge 2,$ it is often of
importance to find an irreducible element passing through a given
set of points. In these circumstances, generality assumptions are
not useful. One is naturally brought to consider sets of points
that are able to break all elements of the given linear system.
The first steps in the treatment of this phenomenon were conducted
in \cite{peculiar1} and \cite{DeF}, where the notions of {\it bad
point} and {\it bad locus} were introduced and studied. Let $L$ be
an ample line bundle on $X,$ spanned by $V \subseteq H^0(X,L).$ A
point $x \in X$ is {\it bad} for the linear system $|V|$ if all
elements of $|V|$ containing $x$ are reducible or non-reduced. The
existence of a bad point for an ample and free linear system
 is shown to be exclusively a two dimensional phenomenon, while bad points do not
occur for very ample linear systems.

The notion of bad point on a projective $n$-fold $X$ can be
generalized in different directions. One can view a single bad
point as a reduced zero-scheme of length one and therefore
generalize the notion to zero-schemes of any length. On the other
hand, recalling that bad points occur only on surfaces, one could
view a bad point as a linear space of codimension two, i.e. a
sub-manifold $\Lambda$ of $X,$  isomorphic to $\Pin{n-2},$ such
that $\restrict{L}{\Lambda}= \oofp{n-2}{1}$.

In this work the more general notions of {\it higher order bad locus,}
 {\it bad zero-scheme,} and {\it bad linear space} are introduced.

A zero-scheme $\xi$ is {\it bad} for the linear system $|V|$ if
all elements of $|V|$ containing (scheme-theoretically) $\xi$ are
reducible or non-reduced. The minimum length of a bad zero-scheme
for the pair $(X,V)$ is introduced as a numerical character
denoted by $b = b(X,V),$ see section \ref{hobl} for details.
Similarly, $b_0$ denotes the minimum length of a bad, reduced
zero-scheme.

The main goal of the first sections of this work is to investigate
relationships among $b_0, b,$ the dimension of $X,$ and the
positivity of $L.$

A crucial point is whether a bad zero scheme of minimal length
imposes independent conditions on $|V|.$ To answer this question
positively, one would need, for any zero-scheme $\xi$ not imposing
independent conditions, to find a subscheme $\eta \subset \xi,$
imposing to $|V|$ linearly independent conditions that are
equivalent to the ones imposed by $\xi.$ While this happens in
several instances, e.g. when $\xi$ is reduced, this fact seems
doubtful in general. To overcome this problem, avoiding
duplication of statements, a notion of {\it suitable} pairs $(V,
\xi)$ is introduced, see section \ref{notation} for details.

In \cite[Theorem 2, i)]{peculiar1} it is shown that if $b(=b_0)=1$
for an ample and free linear system, then $\dim{X} = 2.$ Theorem
\ref{dimbound}, under the assumption that there exists a suitable
pair $(V,\xi),$ gives the bound $\dim{X} \le b+1,$ generalizing
the result above. On the other hand, the corresponding inequality
with $b_0$ instead of $b$ holds with no further assumption, see
Remark \ref{dimbound0}.

Bringing the positivity of $L$ into the picture, one can assume
that $L$ is $k$-very ample, i.\ e., every zero-scheme of length
$k+1$ imposes independent conditions on sections of $L.$ Then,
Theorem \ref{betterbound} gives a stronger bound if $k\ge2.$

In \cite[Theorem 1.7.9]{BESO} a characterization of the case $b = 2$ for
$L$ very ample ($1$-very ample) is given.
In particular $X$ must be a surface and  bad zero-schemes of length $2$ are
contained in a line $\ell.$
Proposition  \ref{b0=k+1} generalizes this results under the assumption that $L$ is $k$-very
ample and $b$ is realized by a reduced zero-scheme of length $k+1.$
In this case $X$ must also be a surface and the bad, reduced zero-schemes of length
$k+1$ are contained in a rational normal curve of degree $k.$

A similar characterization for $k = 2$ and $b=3,$ where $b$ may be
realized by a non-reduced zero-scheme, is given in Proposition
\ref{b=3}.

Sharpness of all bounds is illustrated by a series of examples.

Bad linear spaces are discussed in the final section. It turns out
that they must necessarily have codimension two and that they are
inherited by hyperplane sections. These two facts are combined to
show that bad linear spaces do not occur at all for very ample
linear systems.

The authors are very grateful to the referee for an extremely careful reading of  this work, and in particular for pointing out an inaccuracy in the proofs of earlier versions of Proposition \ref{lowerbound} and \ref{b0=3}.
The first author would also like to thank Lawrence Ein for a useful conversation on issues related to the notion of suitable pairs  introduced in this work.

The authors acknowledge financial support from
G\"oran Gustafsson Grant 2003-2004 for visiting scientists at K.T.H., the
Vetenskapsr{\aa}det 2004-2008 grant, DePaul University, Universit\`a degli Studi di Milano,
and the MIUR of the Italian
Government in the framework of the PRIN Geometry on
Algebraic Varieties, Cofin 2002 and 2004.

\section{Notation and Background}
\label{notation}
 Throughout this article $X$ denotes  a smooth,
connected, projective variety of dimension $n,$ $n$-fold for
short, defined over the complex field $\comp.$ Its structure sheaf
is denoted  by ${\mathcal O}_X$ and the canonical sheaf of
holomorphic $n$-forms on $X$ is denoted by $K_X.$  Cartier
divisors, their associated line bundles and the invertible sheaves
of their holomorphic sections are used with no distinction. Mostly
additive notation is used for their group.

 Let $S^tX$ be the $t$-th symmetric power of $X$ and  $X^{[t]}$
be the Hilbert scheme of zero-subschemes of $X$ of length $t.$
 Let $\red{t}$ be the stratum of reduced zero-subschemes of length $t$.
 We denote by $X^{[t]}_{(1+r,1,1,\dots,1)},$ for $0\le r\le {\rm min}\{t -1,n\},$ the set of zero
 subschemes $\xi$ of length $t$ such that $\supp{\xi} = \{x_1,x_2,\dots, x_{t-r}\},$
  and $\iofo{\xi} = \mathfrak{a}\cdot\mathfrak{m}_2\cdot \ \dots \ \cdot \mathfrak{m}_{t-r}$
 where $\mathfrak{m}_i$ is the maximal ideal of $\calo_{X,x_i}$ and $\mathfrak{a} = (u_iu_j, u_{r+1},\dots,u_n\ |\ 1\le i\le j\le r),$   $u_1,\dots
 u_n$ denoting local coordinates at $x_1.$

Given a zero-scheme  $\xi \in \red{t},$ we sometimes identify $\xi$
and its support $\supp{\xi};$ for example we write $x \in
\xi$ to mean $x \in {\rm Supp} (\xi).$
 For any coherent sheaf ${\mathcal F}$ on $X$,  $h^i(X,{\mathcal F})$ is the complex dimension of
$H^i(X,{\mathcal F}).$ When the ambient variety is understood, we often write $H^i({\mathcal F})$
and $h^i({\mathcal F})$ respectively for $H^i(X,{\mathcal F})$ and $h^i(X,{\mathcal F}).$ Let $L$ be a  line bundle on $X.$ If $L$
is ample, the pair $\xel$ is called a {\it polarized variety}. For
a subspace $V \subseteq H^0 \xel$ the following notations are
used:\\
$|V|,$ the linear system associated with  $V;$\\
$|V \tensor \iofo{Z}|,$ with a slight abuse of notation, the linear system of divisors
 in $|V|$ which contain, scheme-theoretically, the subscheme $Z$ of
 $X;$\\
$\bslocus{V},$ the base locus of the linear system $|V|$;\\
$\phiV,$ the rational map given by $|V|.$ \\
If $V = H^0(L)$ we write $L$ instead of $V$ in all of the above.
Let $L$ be a line bundle generated by its global sections. When
the linear span of $\phiL(\xi)$ is a $\Pin{k}$ for every
zero-scheme $\xi \in X^{[k+1]}$ we say that $L$ is $k$-very ample.
Equivalently:
\begin{defn}
Let $k$ be a non-negative integer. A line bundle $L$ on $X$ is $k$-very ample if the restriction map $H^0(X,L) \to H^0(L \tensor \oofo{\xi})$ is surjective for every zero-scheme $\xi \in \X{k+1}.$
\end{defn}

%
%

The second Bertini theorem, see for example \cite{kl2}, is the
main tool to handle linear systems whose elements are all
reducible. The following remark on the dimension of the base locus of such linear systems follows easily from that theorem and will be useful to us.
\begin{rem}\label{pencil}
Let $L$ be a line bundle on a smooth variety $X,$ spanned by a
subspace $V \subseteq H^0(L).$ Assume $Z$ is a nontrivial
subscheme of $X$ such that $|V \otimes \cI_{Z}|$ does not have a
fixed component and it is composed with a pencil (i.e. the
rational map associated with $|V \otimes \cI_{Z}|$ has
one-dimensional image). Two generic fibers of the pencil are
divisors on $X$ and their intersection lies in the base locus
${\rm Bs}(|V\otimes\cI_{Z}|)$ of $|V \otimes \cI_{Z}|.$ This
implies that $\dim({\rm Bs}(|V\otimes\cI_{Z}|))=n-2$.
\end{rem}
%
%
For ample linear systems we also have the following observation on the dimension of the base locus.
\begin{rem} \label{smallBaseLocus}Let $L$ be an ample line bundle on a smooth variety $X, $
spanned by a subspace $V \subseteq H^0(X,L).$ Assume $\xi$ is a
zero-scheme on $X$ which imposes $k$ independent conditions on
$|V|$, i.e. $\dim(|V\otimes\cI_{\xi}|)=\dim(|V|)-k$. Consider the
cohomology sequence $0\to H^0(L \tensor \iofo{\xi}) \to H^0(L)
\stackrel{r_{\xi}}{\rightarrow} H^0(L\tensor \calo_{\xi}),$ and
let $W = r_{\xi}(V),$ where $\dim{W} = k.$ Choose a basis for $V,$
taking into account the decomposition $V \simeq {\rm Ker}(r_{\xi})
\oplus W$ to describe the map $\phiV.$ Then one sees that $
\phiV({\rm Bs}(|V\otimes\cI_{\xi}|))\subseteq\Proj {W}=\Pin{k-1}$.
As $L$ is ample, $\phiV$ does not contract any positive
dimensional subvariety and thus $\dim({\rm
Bs}(|V\otimes\cI_{\xi}|))\leq k-1.$
\end{rem}
The fact that, given a zero-scheme of length $b$, one can always
find a zero-subscheme of any length $a\leq b$ is going to be
important. For clarity in the exposition we report below a self
contained proof.
\begin{lem}\label{subscheme}
Let $K$ be an algebraically closed field. Let $\xi$ be a
zero-scheme over $K$ of length $b>0$. Let $0<a< b,$ then there is a
zero-subscheme $\xi' \subset\xi$ of length $a.$
 \end{lem}
\begin{pf}
Let $\xi=Spec(A)$. As $\dim(A)=0$, $A$ is an Artinian ring
and thus it is the product of local Artinian rings.
 Let  $A=B_1\times \ \dots \  \times B_k$, where $\supp{B_i}=x_i$ and
 thus $\supp{\xi}=\{x_1,\dots, x_k\}.$  Assume
 $\length{Spec(B_i)}=b_i$, i.e. $b=b_1+\ \dots \ +b_k.$
 If $\xi$ is a reduced subscheme, i.e. $b_1=\ \dots \ =b_k=1$,
 then a subscheme $\xi'$ as desired is given by the
 ring $B_1\times \ \dots \  \times B_a.$

 Assume that $b_i\geq 2$ for some $i$. To prove the assertion of the Lemma it is enough  to show that for every $i$ we can find a subscheme $\eta\subset Spec(B_i)$ with $\length{\eta}=b_i-1.$
 Let $\mathfrak{m_i}$ be the maximal ideal of the local ring $B_i$ and let $h$ be the smallest integer so that $\mathfrak{m_i}^h=0$. Notice that $1<h\leq b_i.$
 Because $\length{Spec(B_i)}=b_i$, it must be $\mathfrak{m_i}^{b_i}=0.$ Otherwise $B_i$ would have a filtration of length $\geq b_i+1:$
 $$ \mathfrak{m_i}^{b_i}\subset \mathfrak{m_i}^{b_i-1} \subset \ \dots \  \subset \mathfrak{m_i}\subset B_i,$$
 which would contradict the assumption that $\dim_K (B_i)=b_i.$
 Let $x\in \mathfrak{m_i}^{h-1}$. The ideal $(x)$ of $B_i,$ generated by $x,$ is a one-dimensional $K$ vector space.
 Consider the surjective map $\phi: B_i\to (x)$ defined by sending $1$ to $x$. Since $\phi(\mathfrak{m_i})=0$ the map $\phi$ factors through $B_i/\mathfrak{m_i}=K$ (because $K$ is algebraically closed) and thus $K$ maps surjectively onto $(x)$.

 This implies that the quotient $B_i/(x)$ defines a subscheme $$\eta = Spec(B_i/(x))\subset Spec(B_i),$$ where
 $$\length{\eta}=\dim_K(B_i/(x))=\dim_K(B_i)-\dim_K(x)=b_i-1.$$\end{pf}

The following two Lemmata deal with changes in positivity of a
$k$-very ample line bundle when blowing-up at reduced
zero-schemes. A detailed study can be found in \cite[4]{DiR5}.

\begin{lem}\label{spanned}
Let  $L$ be a $k$-very ample line bundle on a projective $n$-fold
$X,$ with $k \ge 0.$ Let $\pi:\tilde{X}\to X$ be the blow-up of
$X$ at $\xi\in\red{t}$, for any $0\le t \le k,$ with exceptional
divisors $E_1,\dots,E_t.$ Then the line bundle $\cL=\pi^*(L)-E_1-\
\dots \ -E_t$ is globally generated on $\tilde{X}$.
\end{lem}
\begin{pf}
Let $\cI_{\xi}=\mathfrak{m}_{x_1}\cdot \mathfrak{m}_{x_2}\cdot \
\dots \ \cdot \mathfrak{m}_{x_t}$ be the ideal defining the
reduced zero-scheme $\xi$.  For every point
$y\in(\tilde{X}\setminus\cup E_i)$ consider the zero-scheme of
length $t+1$ defined by the ideal $\cI_{\xi'}=\cI_{\xi}\cdot
\mathfrak{m}_{\pi(y)}.$ It is $\supp{\xi'}=\supp{\xi}\cup \pi(y).$
Because $L$ is $k$-very ample, there is a section of $L$
vanishing at $\supp{\xi}$ and not vanishing at $\pi(y)$. This
implies that there is a section of $\cL$ not vanishing at $y$.

If $y\in E_i$ for some $i$, it corresponds to a tangent direction
$\tau$ to $X$ at $\pi(E_i)=x_i.$  For simplicity let us fix $i=1$.
Choose local coordinates $\{u_1,\dots,u_n\}$ around the point
$x_1$ and assume that the tangent direction corresponds to the
coordinate $u_1$.  The zero-scheme $\xi'\in
X^{[t+1]}_{(2,1,\dots,1)}$ defined by the ideal
$\cI_{\xi'}=(u_1^2, u_2,\dots,u_n)\cdot \mathfrak{m}_{x_2}\cdot \
\dots \  \cdot \mathfrak{m}_{x_t}$ has length $t+1\le k+1,$ and
$\supp{\xi'}=\supp{\xi}$. Because $L$ is $k$-very ample, the map
$$H^0(X,L)\to H^0(L\otimes \calo_{\xi'})$$ is onto.
Therefore there is a section $s\in H^0(X,L)$ which vanishes  at
$\supp{\xi}$ and such that $ds(\tau)\ne 0.$ Let $D = (s)_0$ and
let $s' \in H^0(\cL)$ be the section corresponding to
$\pi^*(D)-E_1-\ \dots \ -E_k.$ Thus $s'(y)\neq 0$.
\end{pf}

\begin{lem}\label{veryample}\cite[4.1]{DiR5}
Let $L$ be a $k$-very ample line bundle on a projective manifold
$X,$ with $k \ge 1.$ Let $\pi:\tilde{X}\to X$ be the blow-up of
$X$ at $\xi\in\red{t}$, for any $t \le k-1,$ with exceptional
divisors $E_1,\dots,E_t.$ Then the line bundle $\cL=\pi^*(L)-E_1-\
\dots \ -E_t$ is very ample.
\end{lem}
\begin{rem}
In the same context as Lemma \ref{veryample}, $\tilde{X}$ can be
naturally identified with a closed subscheme of $X^{[t+1]}$ and,
under this identification, $\cL$ corresponds to a very ample line
bundle of the subscheme. Moreover, if $\dim{X} = 2,$ then such
line bundle extends to very ample line bundle on the whole
$X^{[t+1]}.$ Indeed, let $\xi\in \red{t}$ and for all $l\geq t$
let
$$X^{[l]}_{\xi}=\{ \eta\in X^{[l]}\text{ such that } \supp{\xi}\subseteq \supp{\eta}\}$$
It is $X^{[t]}_{\xi}=\{\xi\}$ and $X^{[t+1]}_{\xi}\cong
\tilde{X}$. Let $s_{\xi}: X\to S^{t+1}X$ be the map defined by
$s_{\xi}(x)=(x,\supp{\xi})$, and $\rho: X^  {[t+1]}\to S^{t+1}X$
be the Hilbert-Chow morphism. It is:
 $$\rho^{-1}(s_{\xi}(X))\cong X^{[t+1]}_{\xi}.$$
Thus the following diagram is commutative:
$$
\begin{diagram}
\node{\tilde{X}}\arrow{e,=}\arrow{s}^{\pi}\node{X^{[t+1]}_{\xi}}\arrow{e}\arrow{s}\node{ X^{[t+1]}}\arrow{s,^f}\\
\node{X}\arrow{e,=}\node{X}\arrow{e}\node{S^{t+1}X}
\end{diagram}
$$

By Lemma \ref{veryample} the line bundle $\cL=\pi^*(L)-E_1-\ \dots
\ -E_t$ defines an embedding of $X^{[t+1]}_{\xi}$.

In the case $\dim(X)=2$ the line bundle $\cL$ can be equivalently
described as follows. Let $\zeta_{t+1}\subset X^{[t+1]}\times X$
be the universal family with projection maps $p_1:\zeta_t \to
X^{[t+1]}$ and $p_2:\zeta_{t+1} \to X$. Consider the rank $(t+1)$
vector bundle $L^{[t+1]}={p_1}_*(p_2^*L)$. Then
$$\cL=\det(L^{[t+1]})|_{X^{[t+1]}_{\xi}}.$$
In other words $\det(L^{[t+1]})$ extends $\cL.$ Note that because
$\dim(X)=2$ and $L$ is $(t+1)$-very ample then the line bundle
$\det(L^{[t+1]})$ is very ample, \cite{cat-go,le}.
\end{rem}
%
%
Let $L$ be a line bundle on a smooth complex projective variety
$X$ of dimension $n \geq 2$ and let $\xi$ be a zero-dimensional
subscheme of $X$. Let $\iofo{\xi}$ be the ideal sheaf of $\xi$ and
consider the exact sequence
$$0 \to L \otimes \iofo{\xi} \to
L \to L \otimes \calo_{\xi} \to 0.$$ Consider the induced
homomorphism $H^0(X,L) \to  H^0(L \otimes \calo_{\xi}).$ We denote
by $r_{\xi}$ its restriction to $V.$

\begin{defn}
\label{suitable}
 Let $X,$ $L,$ $V,$ $\xi,$ and $r_{\xi}$ be as
above. The pair $(V,\xi)$ is {\it suitable} if, whenever $r_{\xi}$
is not surjective, there exists a subscheme $\eta \subset \xi,$
such that $r_{\eta}$ is surjective and $r_{\eta}(V) =
\rho(r_{\xi}(V)),$ where $\rho: H^0(L \tensor \calo_{\xi}) \to
H^0(L \tensor \calo_{\eta})$ is the obvious restriction
homomorphism. This is equivalent to requiring that $\eta$ imposes
linearly independent conditions on $V$ and $|V \tensor
\iofo{\eta}| = |V \tensor \iofo{\xi}|.$
\end{defn}

\begin{lem}
\label{suitableconditions} Let $X,$ $L,$ $V,$ be as in the
definition above. Let $\xi \in X^{[t]}.$ Then $(V,\xi)$ is
suitable in each of the following cases:
\begin{enumerate}
\item[1)] $\xi$ reduced; \item[2)] $\xi \in X^{[t]}_{(1+r,1,1,
\dots, 1)}$; \item[3)] $V = H^0(X,L),$ $L$ is $(k-2)-$very ample
and $t = \length{\xi}\le k.$
\end{enumerate}
\end{lem}
\begin{pf}
Assume that  $r_{\xi}$ is not surjective, i.e. $\xi \in X^{[t]}$
does not impose $t$ independent conditions. To prove {\it 1)}, let
$\supp{\xi} = \{x_1,\dots,x_t\}$ and let
$\cI_{\xi}=\mathfrak{m}_1\cdot \ \dots \ \cdot\mathfrak{m}_t$
where $\mathfrak{m}_i$ is the maximal ideal of $\calo_{X,x_i}.$
Consider the vector subspace
$${\rm Im}(\,r_{\xi})\subset H^0(L\tensor\oofo{\xi})\cong \oplus_1^{t} \oofo{X}/\mathfrak{m}_i\cong \comp^{t}.$$
After extending a basis of ${\rm Im}(\,r_{\xi})$ to
$H^0(L\tensor\oofo{\xi})$ one can assume that ${\rm
Im}(\,r_{\xi})\cong \oplus_1^{s} \oofo{X}/\mathfrak{m}_i\cong
\comp^{s}$, where $s < t.$ Then the reduced zero-subscheme
$\eta=\{x_1,\dots,x_s\}$ imposes independent conditions on $V$
and it is what we need.

To prove {\it 2)} let $\supp{\xi} = \{x_1,\dots,x_{t-r}\}$ and let
$u_1,\dots,
 u_n$ denote local coordinates at $x_1.$
 Recall that $\iofo{\xi} = \mathfrak{a}\cdot\mathfrak{m}_2\cdot \ \dots \
\cdot \mathfrak{m}_{t-r}$
 where $\mathfrak{m}_i$ is the maximal ideal of $\calo_{X,x_i}$
 and $\mathfrak{a} = (u_iu_j,u_{r+1},\dots,u_n\ |\ 1\le i \le j\le r).$
  It is $$H^0(L \tensor
 \calo_{\xi})\cong \comp^t,$$ where the first summand on the right hand side, which is isomorphic to $\calo_X/\mathfrak{m}_1,$
 is contained in $\rm{Im}(\,r_{\xi})$ as $V$
 spans $L.$
Notice that the restriction map $\rho:H^0(X,L)\to H^0(L \tensor
 \calo_{\xi})$ acts on global sections $s \in H^0(X,L)$ as follows: $\rho(s) =(s(x_1),\ddui{s}{1}(x_1), \dots, \ddui{s}{r}(x_1), s(x_2), \dots, s(x_{t-r})).$ Because only first derivatives at $x_1$ appear,
one can proceed as in the proof of {\it 1)} by completing the
generator of the first summand above to a basis of
$\rm{Im}(\,r_{\xi}).$

To prove 3), note that as $L$ is $(k-2)$ very ample, then
necessarily $t = k.$ By Lemma \ref{subscheme} there exists a
subscheme $\eta \subset \xi$ of length $k-1.$ Consider the
following commutative diagram:
$$
\begin{diagram}
\node{H^0(X,L)}\arrow{e,t}{r_{\xi}}\arrow{s,=}\node{H^0(X, L \otimes \calo_{\xi})\cong \comp^k}\arrow{s,r}{\rho}\\
\node{H^0(X,L)}\arrow{e,t}{r_{\eta}}\node{H^0(X, L \otimes
\calo_{\eta})\cong \comp^{k-1}.}
\end{diagram}
$$
As $L$ is $(k-2)$-very ample, $r_{\eta}$ is surjective. Therefore
$H^0(X, L \otimes \calo_{\eta}) \cong \rm{Im}(\,r_{\xi}),$ and
$\eta$ is the required subscheme.
\end{pf}

The above Lemma is mostly intended as a useful tool to enhance
readability of a number of results expressed in terms of suitable
pairs contained in sections \ref{hobl} and \ref{blhob}.
%
%
\section{Higher Order Bad Loci}
\label{hobl}
The definition of bad locus, introduced in \cite{peculiar1} and
further studied in \cite{DeF}, can be fairly naturally generalized
to subsets of $X^{[t]}$ as follows.
\begin{defn}
\label{thbadlocus} Let $X$ be a complex, non singular, projective
variety. Let $L$ be a  line bundle  on $X$ spanned by a subspace
$V \subseteq H^0(X,L).$
\begin{enumerate}
\item  The $t$-th bad locus of $(X,V),$ for $t\geq 1$ is:
\begin{alignat}{1}\cB_t(X,V)=\{\xi\in X^{[t]}\ |\ &|V\otimes\cI_{\xi}|\neq\emptyset \text{ and  }\notag
\\ &\forall D\in |V\otimes\cI_{\xi}|,\, D\text{ is reducible or non-reduced}\}.\notag \end{alignat}
\item The reduced $t$-th bad locus of $(X,V),$ for $t\geq 1$ is:
$$\cB_t^0(X,V)=\cB_t(X,V)\cap \red{t}.$$
\end{enumerate}
We write $\cB_t(X,L)$ and $\cB_t^0(X,L)$ if $V = H^0(X,L).$
\end{defn}
An element $\xi \in \badV{t}$ is called a {\it bad} zero-scheme for the linear system $|V|.$
%
%

 There is a clear relationship among the $\cB_t$'s:
\begin{lem}
If $\cB_t(X,V)\neq\emptyset$ then $\cB_k(X,V)\neq\emptyset$ for every $k\geq t$.
\end{lem}
\begin{pf} Let $\xi\in \cB_t(X,V)$.
 If $\dim| V\otimes\cI_{\xi}|\geq 1$ let $x$ be any point in $ X\setminus \supp{\xi}.$
Otherwise let $D$ be the unique element in $|V\otimes\cI_{\xi}|$ and let $x$ be any point in
$D\setminus \supp{\xi}.$ Consider the zero-scheme $\xi'$, of length $t+1,$
obtained by adding the reduced point $x$ to $\xi.$ It is $\supp{\xi'}=\supp{\xi}\cup\{x\}$,
 ${\mathcal O}_{\xi',y}={\mathcal O}_{\xi,y}$ for every $y\in \supp{\xi}$ and ${\mathcal
 O}_{\xi',x}=\frac{\calo_X}{\mathfrak{m}_x}.$
Because $|V\otimes\cI_{\xi'}|\subseteq |V\otimes\cI_{\xi}|$ all
the divisors in $|V\otimes\cI_{\xi'}|$ are reducible and $\xi'\in
\cB_{t+1}(X,V)$.\end{pf}

The above Lemma suggests the following
definitions.
\begin{defn}
Let $(X,V)$ be as above. The $b$-index of the pair $(X,V)$ is:
\[ b(X,V)=\left\{\begin{array}{ll}
\infty&\text{ if } \cB_t(X,V)=\emptyset \text{ for every }t\geq 1\\
\min\{t\ |\ \cB_t(X,V)\neq\emptyset\}& \text{ otherwise.}
\end{array}\right.\]
The reduced $b$-index of the pair $(X,V)$ is:
\[ b_0(X,V)=\left\{\begin{array}{ll}
\infty&\text{ if } \badO{t}(X,V)=\emptyset \text{ for every }t\geq 1\\
\min\{t\ |\ \badO{t}(X,V)\neq\emptyset\}& \text{ otherwise.}
\end{array}\right.\]
 We write $b$
and $b_0,$ respectively, for $b(X,V)$ and $b_0(X,V),$ when the
pair $(X,V)$ is clear from the context.
\end{defn}

\begin{rem}
\label{bleb0} It follows immediately from the above definition
that $b(X,L) \le b_0(X,L).$ Moreover, $\cB_b(X,V) \cap \red{b}\neq
\emptyset$ if and only if  $b(X,V) =  b_0(X,V).$ Notice also that
if there are no suitable pairs $(V,\xi),$ then $b < b_0.$
\end{rem}

\begin{rem}\label{picz}If $|V|$ contains a reducible element $D,$ then $b(X,V)
< \infty.$ Indeed, let $A$ be an irreducible component of $D.$ A
zero-scheme $\xi \subset A$ can be constructed with $r =
\length{\xi}
> \dim{\restrict{|V|}{A}}$ and  sufficiently
general to have $|V \tensor \iofo{\xi}| = A + |V - A|.$ Then $\xi
\in \badV{r},$ hence $b(X,V) \le r.$

Suppose $\Picof{X} = \mathbb{Z}[L]$ where $L$ is ample and spanned
by $V.$ Then $b(X,V) = \infty.$ Indeed $|V|$ cannot contain any
reducible element. If $A + B$ were such an element it would be $A
= aL$ and $B = bL$ for $a,b \ge 1.$ This would give $L = (a+b)L$
which is a contradiction. Recall that Barth-Larsen's Theorem can
provide plenty of such examples.
\end{rem}

\begin{rem} Let $W \subseteq V \subseteq H^0(X,L),$ be two subspaces
which generate $L.$ Then clearly $\cB_t(X,V) \subseteq \cB_t(X,W)$
and thus $b(X,W)\le b(X,V)$ and $b_0(X,W)\le b_0(X,V).$ This is illustrated in the following example.
\end{rem}

\begin{exmp}
\label{p202WV} Let $(X,L)=(\Pin{2}, \oofp{2}{2})$ and $\xi \in
X^{[2]}.$  As $|L \tensor \iofo{\xi}|$ contains always irreducible
conics, it is $b (X,L) \ge 3.$ Let now $\eta$ consist of three
distinct points on a line. As all elements of $|L \tensor
\iofo{\eta}|$ are reducible, it is $b(X,L)= b_0(X,L)= 3.$ Now let
$x_0, x_1, x_2$ be homogeneous coordinates on $\Pin{2}$ and
consider the following vector subspace of $H^0(X,L)$: $U:=\langle
x_0^2+x_1^2+x_2^2, x_0x_1, x_0x_2, x_1x_2 \rangle$. Note that $U$
spans $L$ and $\varphi_U:\Pin{2} \to \Pin{3}$ is a birational
morphism whose image, $\Sigma,$ is Steiner's roman surface.
Theorem 1.1 in \cite{DeF}, since $\varphi_U(X)$ is neither
$\Pin{2}$ nor a cone, implies $\cB_1(X,U) = \emptyset.$ Notice
that $\varphi_U$ maps all three points $e_0=(1:0:0), e_1=(0:1:0),
e_2=(0:0:1)$ to $(1:0:0:0)$, the triple point of $\Sigma$. In
other words,
$$|U-e_0|=|U-e_1|=|U-e_2|.$$ Now, let $\xi$ be the zero-scheme
consisting of $e_i$ and another point $p$, possibly infinitely
near, lying on the line $\langle e_i, e_j \rangle$. Then every
conic in $|U \otimes \cI_{\xi}|$ is reducible, containing $e_i,
e_j, p$, hence the line $\langle e_i, e_j \rangle$. This shows
that $\xi \in \cB_2(X,U)$  and therefore $b(X,U)= b_0(X,U)=2 < b(X,L)$. Finally, let $W$
be either the vector subspace $\langle x_0^2, x_1^2,
f(x_0,x_1,x_2) \rangle$, where $f$ is a general form of degree
$2$, or $\langle x_0^2, x_1^2, x_2^2 \rangle$. Then the pairs
$(X,W)$ correspond to Example 7.2 (i) and Example 7.3 (jjj) in
\cite{DeF}, respectively. Recalling that $\cB(X,W)= \{e_2\}$
or $\{e_0,e_1,e_2\}$ respectively, in both cases we have $b(X,W)=
b_0(X,W) = 1 < b(X,U)< b(X,L)$.
\end{exmp}

The following Proposition generalizes some points of \cite[Theorem
2]{peculiar1}.
\begin{prop}
\label{multiplicity}
 Let $X$ be a smooth  $n$-dimensional variety. Let $L$ be an ample
line bundle on $X,$ spanned by a subspace $V \subseteq H^0(X,L).$
Let $\xi \in \bad{b}(X,V),$ and suppose $|V \otimes \cI_{\xi}|$
has finite base locus. Then
\begin{enumerate}
\item[i)] $n=2;$
\item[ii)]there is an ample line bundle $A$ on $X$ with $h^0(A)\ge
2$ such that  every $D \in |V \otimes \cI_{\xi}|$ is of the form
$D = A_{b_1} + \ \dots \  + A_{b_r},$ for some $r\ge2,$ with $A_{b_i}$ varying in a
rational pencil $B\subseteq|A|;$
\item[iii)] with $r$ as in $ii),$ for all $x \in \bslocus{V \otimes \cI_{\xi}},$  $x$ is a point of multiplicity $r \ge 2$ for all $D
\in |V \otimes \cI_{\xi}|.$ In particular this is true for all $x \in \supp{\xi},$ i.e. $\iofo{\xi,x} \subseteq
\mathfrak{m}_x^r,$ where $r \ge 2;$
 \item[iv)] $\supp{\xi}\subseteq A_{b_i} \cap
A_{b_j}=\bslocus{V \otimes \cI_{\xi}},$  for all distinct $A_{b_i}$ and $A_{b_j}$
appearing in the expression of a general $D \in |V \otimes
\cI_{\xi}|.$
\end{enumerate}
\end{prop}
\begin{pf}
As $|V \otimes  \cI_{\xi}|$ has finite base locus, Bertini's
second theorem and Remark \ref{pencil} give $i)$ and that the
image $C$ of the rational map $\varphi_{V\otimes \cI{\xi}}$ is
one-dimensional. After resolving its indeterminacies, taking
Stein's factorization we get the
 following diagram:
 $$
 \begin{diagram}
\node{\tilde{X}}\arrow{e,t}{\tilde{\varphi}}\arrow{s,t}{\alpha}\node{C\subset \Pin{}}\\
\node{B}\arrow{ne,t}{\beta}
\end{diagram}
 $$
 where $\tilde{X}$ is a suitable blow-up of $X,$ $B$ is a smooth
 curve, $\alpha$ has connected fibres and $\beta$ is a finite
 morphism. Note that $B \simeq \Pin{1}$ because there is at least
 one exceptional divisor in $\tilde{X},$ mapping surjectively to $B$ via
 $\tilde{\varphi}.$
 Let $r = \deg{\beta} \deg{C}.$ Then every $D \in |V \otimes
 \cI_{\xi}|$ can be written as $D = A_{b_1} + \ \dots \  + A_{b_r},$
 where each $A_{b_i}$ is the image on $X$ of a fibre of $\alpha.$
 Thus, $r\ge 2.$ Notice that all $A_{b_i}$'s are linearly
 equivalent as they vary in the rational pencil $B.$ It follows
 that $h^0(A) \ge 2,$ $L \sim rA$ (linearly equivalent) and $A$ is
 ample. This proves $ii).$
 Let $x\in \bslocus{V \otimes
 \cI_{\xi}}.$ As $ x \in D$ for all $D \in |V \otimes
 \cI_{\xi}|,$ and there are no fixed components, $x$ must belong
 to infinitely many elements  $A_{b}, b \in B,$ hence to all of
 them. Thus $x$ is a point of multiplicity greater or equal to $r$
 for all $D.$ In particular, if $x\in \supp{\xi},$ $\iofo{\xi,x} \subseteq \mathfrak{m}_x^r.$
This proves $iii).$ Moreover $\supp{\xi}\subseteq \bslocus{V \otimes
 \cI_{\xi}} \subseteq
A_{b_i}\cap A_{b_j}$ for any $i,j =1, \dots, r.$ To prove $iv),$
it is then enough to show that for a general $D \in |V \otimes
\iofo{\xi}|$ it is $\bslocus{V \otimes
 \cI_{\xi}} \supseteq A_{b_i}\cap A_{b_j}$ for
any $i,j =1, \dots, r.$ This follows from Bertini's first theorem
because every point $y \in A_{b_i}\cap A_{b_j}$ is a singular
point for $D,$ which is generally chosen.

\end{pf}
The following Lemma shows that, in the case of suitable pairs, the linear span of a bad zero-scheme of minimal length is always of maximal dimension.

\begin{lem}\label{conditions} Let $X$ be a smooth  $n$-dimensional variety. Let $L$ be an ample
line bundle on $X,$ spanned by a subspace $V \subseteq H^0(X,L).$
Assume $b_0(X,V)< \infty.$
\begin{enumerate}
\item[i)]Let $\xi\in \bad{b}(X,V).$ If $(V,\xi)$ is suitable, then
$\xi$ imposes exactly $b$ independent conditions on $|V|;$
\item[ii)] $b(X,V) \le b_0(X,V) \le \dim |V|.$
\end{enumerate}
\end{lem}
\begin{pf}
 To prove $i),$ assume that
$\dim(|V\otimes\cI_{\xi}|)=\dim(|V|)-m>\dim(|V|)-b$. As $(V,\xi)$
is suitable, there exists a zero-scheme $\eta \subset \xi$ with
$\length{\eta} < b,$ such that
$|V\otimes\cI_{\eta}|=|V\otimes\cI_{\xi}|$. This is impossible
because it would imply $\eta \in \bad{t}(X,V)$ with $t < b$. If
$\xi$ in the argument above is reduced then  $b = b_0,$ as noted
in Remark \ref{bleb0}. Moreover in this case  $\eta$ is reduced
and the right side of inequality $ii)$ follows immediately from
$i).$ Remark \ref{bleb0} completes the proof.
\end{pf}

\begin{rem}
\label{conditions0} Notice that replacing $\badV{b}$ with
$\badOV{b_0}$ and $b$ with $b_0,$ the same argument as in the
proof of Lemma \ref{conditions} shows that $i)$ holds for every
$\xi \in \badOV{b_0}.$
\end{rem}
The upper bound in Lemma \ref{conditions} $ii)$ above can be
strict, see example \ref{p202WV} above with $V = H^0(X,L)$
 and example \ref{p3o2} below. Nonetheless it
can be attained in some cases, showing that it is optimal. This
can be seen in example \ref{antonio1}.

\begin{exmp}
\label{p3o2} Let $X=\Pin{3}$ with homogeneous coordinates $x_0,
x_1, x_2, x_3,$ and consider the vector subspace $V$ of
$H^0(\oofp{3}{2})$ generated by the monomials $x_0^2, x_1^2,
x_2^2, x_3^2$. Note that $V$ spans $L=\oofp{3}{2}$ and defines a
morphism $\Pin{3} \to \Pin{3}$ of degree $8$. Let
$\{e_0,e_1,e_2,e_3\}$  be the standard basis for $\comp^4.$ Now,
let $\xi = \{e_0,e_1\}.$ Then $|V \tensor \iofo{\xi}|$ is the
pencil of quadrics generated by $x_2^2$ and $x_3^2$. Every such
quadric is reducible, hence $\xi \in  \badOV{2}$. Note that $|V
\tensor\iofo{\xi}|$ has no fixed component and is composed with
the pencil $|\oofp{3}{1}\tensor\iofo{\ell}|$, $\ell$ being the
linear span of $\xi$. In conclusion we have
$$b(X,V) \le b_0(X,V)= 2 =  \dim{|V|} -1,$$ and thus $b(X,V) = 2,$
as $\cB(X,V) = \emptyset,$ see \cite[Theorem 2]{peculiar1}.
\end{exmp}

\begin{exmp}
\label{antonio1} Let $X$ be a Del Pezzo surface with $K_X^2=2$ and
let $L=-K_X$. It  is well known that $L$ is ample and spanned and
$\phiL:X \to \Pin{2}$ is a double cover, branched along a smooth
plane quartic curve. Let $\xi =\{x_1, x_2\}$ be a zero-scheme consisting of two distinct
points of the ramification divisor such that $<\phiL(x_1),
\phiL(x_2)>$ is a bitangent line to the branch quartic curve.
Then $|L \tensor \iofo{\xi}|$ consists of a single element $G$ having double
points at $x_1$ and $x_2$. In fact $G = \Gamma_1+\Gamma_2$, where
$\Gamma_1$, $\Gamma_2$ are two $(-1)$-curves meeting exactly at
$x_1, x_2$. Thus $b_0(X,L) \leq 2$ and it follows from Corollary 1.3 in
\cite{DeF} that $b_0(X,L)=2=\dim |L|$.
\end{exmp}
One of the major results of \cite{peculiar1} is the fact that the existence of a bad
point forces the variety to be a surface.
The following Theorem generalizes this result to include higher order bad loci
\begin{thm}
\label{dimbound} Let $X$ be a smooth  $n$-dimensional variety. Let
$L$ be an ample line bundle on $X,$ spanned by a subspace $V
\subseteq H^0(X,L).$  Assume there exists $\xi \in
 \badV{b}$, such that $(V,\xi)$ is suitable. Then:
\begin{enumerate}
\item[1)] $n\leq b+1$; \item[2)]If $n=b+1$, for every $\xi\in
\bad{b}(X,V)$ such that $(V,\xi)$ is suitable, the linear system
$|V\otimes\cI_{\xi}|$ has no fixed component and it is composed
with a pencil.
\end{enumerate}
\end{thm}
\begin{pf}
Let $\xi\in \bad{b}(X,V)$. If $b =1,$ {\it 1)} and {\it 2)} follow
respectively from \cite[Theorem 2, i) and ii)]{peculiar1}. So let
$b\ge 2.$ By Lemma \ref{conditions}, $\xi$ imposes exactly $b$
conditions on $|V|$. According to Bertini's second theorem, the
linear system $|V\otimes\cI_{\xi}|$ either has a fixed component
$\Sigma$ or it is composed with a pencil. In the first case,
Remark \ref{smallBaseLocus} gives $\dim(\Sigma)=n-1\leq b-1,$ i.e.
{\it 1)} holds as a strict inequality and thus {\it 2)} is proven.
If $|V\otimes\cI_{\xi}|$ is composed with a pencil, Remarks
\ref{pencil} and \ref{smallBaseLocus} give
$$n-2 = \dim({\rm Bs}(|V\otimes\cI_{\xi}|))\leq b-1,$$
which completes the proof of {\it 1)}.
\end{pf}
\begin{rem}
\label{dimbound0}
 If there are no suitable pairs $(V,\xi),$ which
by Remark \ref{bleb0} implies $b<b_0,$ the same argument used in
the proof of Theorem \ref{dimbound}, replacing $\badV{b}$ with
$\badOV{b_0}$ and Lemma \ref{conditions} with Remark
\ref{conditions0}, gives the same statement as in Theorem
\ref{dimbound} for $b_0$ and for all $\xi \in \badOV{b_0}.$
\end{rem}
Example \ref{p3o2} shows that case {\it 2)} of Theorem
\ref{dimbound} is effective. On the other hand the following
example illustrates the fact that if $n < b+1$ then the linear
system $|V \tensor \iofo{\xi}|$ for all $\xi \in \bad{b} (X,V)$
may have a fixed component.
\begin{exmp}
\label{containsLplane} Let $X$ be an $n$-fold and let $L$ be an
ample line bundle on $X,$ spanned by a subspace $V \subset
H^0(X,L).$ \\ Assume that $X$ contains an $L$-hyperplane, i.e. a
divisor $F \simeq \Pin{n-1}$ such that $L^{n-1}\cdot F = 1.$ Let
$\xi =\{ x_1, x_2,\dots, x_n\}$ be a reduced zero-scheme on $F$
such that the linear span of $\phiV(\xi)$ has dimension $n-1,$
i.e. $\phiV(\xi)$ spans the entire $\phiV(F).$ Then $\xi \in
\badO{n}$ with $|V \otimes \iofo{\xi}|$ having $F$ as fixed
component. It follows that $b (X,V) \le b_0(X,V)\le n.$ On the
other hand Theorem \ref{dimbound} gives $b(X,V) \ge n,$ hence
$b(X,V) = n.$

In particular, if $\xel = \scroll{E}$ is the $n$-dimensional
scroll of a very ample vector bundle $E$ over a smooth curve $C$ then
$b (X,L) = b_0 (X,L) = n.$ Notice that as $(X,L)$ is a scroll, $D\in |L|$ is reducible
if and only if $D$ contains a fibre $F = \Pin{n-1}$ of the scroll.
\end{exmp}

\section{Bad loci and higher order embeddings}
\label{blhob}
 The study of $\cB_t(X,L),$ so tightly connected with
linear systems containing zero-schemes of any length, is very
naturally
 conducted in the context of higher order embeddings. The focus of
 this section is on the properties of bad loci of complete
 linear systems associated with $k$-very ample line bundles.

 Let $X$ be a smooth complex variety and let $L$ be a $k$-very
ample line bundle on $X.$ Assume $\cB_{t}(X,L)$ is not empty for
some $t.$ Recall that, by Lemma \ref{suitableconditions}, for all
$\xi \in X^{[t]},$ $t \le k+2,$ the pair $(H^0(X,L), \xi)$ is suitable.

In \cite{peculiar1} and \cite{DeF}, the case of a spanned
($0$-very ample) line bundle with non empty $\cB_1(X,L)$ was
treated.

When $L$ is very ample ($1$-very ample) it is $b_0\ge b \ge 2,$ see
\cite[Corollary 2]{peculiar1},
and a complete characterization of the case $b = 2$ is given in
 \cite[Theorem 1.7.9]{BESO}.
In this case $n=2$ and, for all $D \in |L|$ containing a bad zero-scheme of length 2, $D=\ell + R$,
 where $\ell$ is a fixed line containing the bad zero-scheme,
$\oof{X}{R}$ is spanned, and $R$ and $\ell$ intersect transversally.
Notice that the same characterization holds for $(X,V),$ $V \subset H^0(X,L),$
where $|V|$ is a very ample linear system.
Notice also that this shows that the natural phenomenon described in
Example \ref{containsLplane} is the only possibility when $L$ is very ample and $n = 2.$
The following Proposition generalizes the lower bound on $b$ in terms of $k$-very ampleness.
\begin{prop}\label{lowerbound}
Let  $L$ be a $k$-very ample line bundle on a projective $n$-fold
$X,$ $n \ge 2,$ with $k \ge 1.$ Then $b(X,L) \geq k+1$ if either
$n\ge 3$ or $n=2$ and there exists $\eta \in \cB_b(X,L)$ with $x
\in \supp{\eta}$ such that $\eta$ is reduced at $x,$ i.e.
$h^0(\calo_{\eta,x}) = 1.$
\end{prop}
\begin{pf} Assume by contradiction that $b \leq k$ and let $\xi \in
\cB_b(X,L)$. Notice that the base scheme of $|L \otimes
\iofo{\xi}|$ is $ \xi.$ Indeed, if such a scheme, $Z,$ strictly
contained $\xi$ as a scheme, then there would be a zero-scheme
$\xi'$ of $X$, containing $\xi$, of length $b+1$ such that $|L
\otimes \iofo{\xi}|=|L \otimes \iofo{\xi'}|$. But this contradicts
the $k$-very ampleness of $L$, since $b+1 \leq k+1$. Then the
assertion follows from Proposition \ref{multiplicity}, $i)$ and
$iii).$
\end{pf}

\begin{rem}
\label{p20k} Equality in the statement of Proposition
\ref{lowerbound} does indeed happen. Let $\xel = (\Pin{2},
\oofp{2}{k}).$ If $k=1$ then $b(X,L) = \infty$ by Remark
\ref{picz}. Assume $k \ge 2.$ $L$ is $k$-very ample. Let $\xi =
\{x_1, \dots , x_{k+1}\}$ be a reduced zero-scheme contained in a
line $\ell.$ Then $\xi \in \badO{k+1}\xel$ and thus Proposition
\ref{lowerbound} gives $b_0 \xel = b \xel = k+1.$
\end{rem}

In line with \cite[1.7.9]{BESO} where the case $b=2$ is fully
described, the following Proposition gives a complete
characterization of the case $b=3.$ The role of the line $\ell$
for $b=2$ is here played by a smooth conic.

\begin{prop}
\label{b=3}
 Let  $L$ be a $k$-very ample line bundle on a
projective $n$-fold $X,$ with $n \ge 2$ and $k \ge 2.$ Assume
$b(X,L) = 3.$ Then  $n = k = 2$ and either
\begin{enumerate}
\item[a)] $(X,L) =(\Pin{2}, \oofp{2}{2})$ or \item [b)] for all
$\xi \in \badL{3},$  $|L \tensor \iofo{\xi}|$ has a fixed
component $\Gamma$ which is a smooth conic containing $\xi.$

\end{enumerate}
\end{prop}
\begin{pf}
We first show that $k=2.$ Let $\xi \in \badL{3}.$ Because
$b(X,L)=3,$ Proposition \ref{lowerbound} gives  $k = 2$ except
possibly when $n=2$ and $\xi$ is supported on a single point. By
contradiction, assume $n=2, k\geq 3,$ and let $\supp{\xi}=\{y\}.$
The same argument as in the proof of Proposition \ref{lowerbound}
shows that $\xi$ is the base scheme  of $|L\otimes \iofo{\xi}|$.
Then, by Proposition \ref{multiplicity}, all $D\in |L\otimes
\iofo{\xi}|$ are of the form $D = A_{b_1} + \ \dots \  + A_{b_r},$
 where $A_{b_i}$  varies in a
rational pencil $B\subseteq |A|,$ and $r=2$ as $b(X,L)=3.$ Notice
that, for a general $D,$ $A_{b_i}$ is smooth at $y.$ Otherwise $y$
would be a point of multiplicity $\ge 4$ for all $D,$
contradicting $b(X,L) = 3.$ Because $h^0(\calo_\xi) =
h^0(\calo_{\xi,y})=3,$ it is $A^2\le 1$ and hence $A^2 = 1,$ $A$
being ample. Indeed if $A^2\ge2$ then  Proposition \ref{multiplicity}, $iv),$ noting that $\supp{\xi} = \bslocus{L \otimes \iofo{\xi}},$  implies that all $A_{b_i}$ have an
assigned tangent at $y.$ Then, locally at $y,$ each
 $A_{b_i}$ has an equation of the form $z_2+F_i(z_1,z_2)=0$ where $F_i$
 has no terms of degree less than $2.$
Therefore $H^0(\calo_{\xi,y})\supseteq <1,z_1,z_2,z_1^2,z_1z_2>,$
which is a contradiction.
 It follows that $L$ is a $3$-very
ample line bundle, with $L^2=(2A)^2=4,$ which is impossible. To
see this, blow-up $X$ at  two general points $x_1,x_2\in X,$ and,
on the new surface $\tilde{X},$ consider the very ample line
bundle $\mathcal{L}=\pi^*(L)-E_1-E_2$ as in Lemma \ref{veryample}.
Note that $\mathcal{L}^2= 2$  while the Picard number of
$\tilde{X}$  is at least $3,$ which is clearly impossible. Thus
$k=2.$ The proof now splits into $2$ cases, according to the
cardinality of $\supp{\xi}.$

{\it Case} 1. Assume first $\xi = \xiP \cup \{x\}$ where
$\length{\xiP} = 2$ and $ x \not \in \supp{\xiP}.$ Let
$\pi:\tilde{X}\to X$ be the blow-up of $X$ at $x,$ with
exceptional divisor $E.$ Then the line bundle $\cL=\pi^*(L)-E$ is
very ample by Lemma \ref{veryample}. Let $\xiPT$ be $\xiP$ pulled
back on $\tilde{X}.$ Note that there is a bijection between $|\cL
\tensor \iofo{\tilde{\xiP}}|$ and $|L \tensor \iofo{\xi}.|$ Hence
$\xiPT \in \cB_2(\tilde{X}, \cL).$ From \cite[1.7.9]{BESO} it
follows that $n=2$ and $|\cL \tensor \xiPT| = \ell + |R|,$ where
$\ell$ is an $\cL$-line containing $\xiPT$ scheme theoretically.
Let $\Gamma = \pi(\ell).$ Because $L$ is $2$-very ample it must be
$L\cdot \Gamma \ge 2,$ with equality holding only if $\Gamma
\simeq \Pin{1}.$ Because $\ell$ and $E$ are $\cL$-lines, it must
be $\nu:=\ell \cdot E \in\{ 0,1\}.$ Therefore it is $1 = \cL \cdot
\ell = (\pi^*(\Gamma)-  \nu E)\cdot (\pi^*(L) - E) = \Gamma \cdot
L - \nu \ge 2-1 = 1$ which implies $\Gamma \cdot L = 2,$ i.e.
$\Gamma$ is a smooth conic. Moreover $\ell \cdot E = 1,$
 hence $x$ belongs to $\Gamma.$ Thus we are in case $b)$.

{\it Case} 2. Assume now that $\supp{\xi} = \{x\}.$ Let
$\pi:\tilde{X}\to X$ be the blow-up of $X$ at $x,$ with
exceptional divisor $E.$ As above the line bundle $\cL=\pi^*(L)-E$
is very ample. Note that $\xi$ defines a length $2$ zero-scheme
$\eta$ of $\tilde{X}$ supported on  $E.$ Again
 there is a bijection between $|\cL \tensor \iofo{\eta}|$ and $|L \tensor
\iofo{\xi}|,$ and thus  $\eta \in \cB_2(\tilde{X}, \cL).$ As
before we conclude that $n=2$ and $|\cL \tensor \iofo{\eta}| =
\ell + |R|,$ where $\ell$ is an $\cL$-line containing $\eta$
scheme theoretically. If $\ell\neq E$ then as above $\Gamma =
\pi(\ell)$ is a smooth conic containing $x$ (hence $\xi$) and it
is a fixed component for $|L \tensor \iofo{\xi}|.$ We are again in
case $b)$.
 Let now $\ell=E$.  In view of the
quoted bijection $|L \tensor \iofo{\xi}|$ corresponds to $|R|$ and
thus all elements of $|R|$ are reducible. According to
\cite[1.7.9]{BESO}, $|R|$ is base-point free, hence its generic
element is smooth. By Bertini's second theorem $|R|$ is composed
with a pencil and thus $R^2=0.$  The pencil $B$ is rational
because $E$ is transverse to all fibers. From $1=\cL\cdot
\ell=\ell^2+\ell\cdot R=-1+\ell\cdot R$ we get that $\ell\cdot
R=2,$ hence $R=A_b + A_{b^\prime}, b,b^\prime \in B, \cL \cdot A_b
= \cL \cdot A_{b^\prime} = 1.$  This shows that $\tilde{X}$ is
fibred over $\Pin{1}$ by $\cL$-lines. Moreover
$\cL^2=\cL\cdot(\ell+R)=3,$ and then it follows that
$(\tilde{X},\cL)$ is a rational cubic scroll. Thus
$(X,L)=(\Pin{2}, \oofp{2}{a}),$ for some $a\ge 2.$ Let $f$ be a
fiber of the cubic scroll. Then $\pi^*(\oofp{2}{1})=E+f.$ On the
other hand $R=\cL-E$, which implies that $2E+af=\pi^*L=a(E+f)$ and
thus $a=2.$ This gives case $a)$.
\end{pf}

The following example illustrates case $b)$ of Proposition
\ref{b=3}.
\begin{exmp}
\label{antonio1extended} Let $X$ be a non-minimal Del Pezzo
surface with $K_X^2 \geq 2$ and let $L=-2K_X$. Then $L$ is
$2$-very ample, see \cite{DiR2}. Let $E \subset X$ be a
$(-1)$-curve, and let $\xi \in X^{[3]},$ supported on $E$. Note
that $|L \otimes \cI_{\xi}|$ is nonempty, since $H^0(X,L)=
1+3K_X^2 \geq 7$. Moreover, it is clear that $E$ is a fixed
component of $|L \otimes \cI_{\xi}|$, since $LE=2$. This shows
that $\xi \in \cB_3(X,L)$. It thus follows that $b(X,L) \leq 3$.
Note however that it cannot be $b(X,L)=2$, since $(X,L)$ does not
contain lines. Therefore $b(X,L)=3$. Note also that the fixed
component $E$ of $|L \otimes \cI_{\xi}|$ is a rational normal
curve of degree $2$ in the embedding given by $L$.
\end{exmp}

The following example shows that when $L$ is very ample, but not
$2$-very ample, the fixed component of $|L \tensor \iofo{\xi}|$
can be singular.

\begin{exmp}
\label{antonio2} Let $X$ be a Del Pezzo surface with $K_X^2=1$ and
let $L=-3K_X$. Then $L$ is very ample. However $L$ is not $2$-very
ample, see \cite{DiR2}. Recall that if $X$ is general in moduli
then $|-K_X|$ is a pencil containing 12 irreducible elements
having a double point. Let $\Gamma$ be such an element and let
$x_1$ be its singular point. Let $x_2, x_3$ be two other distinct
points of $\Gamma$ and consider the reduced zero-scheme $\xi$ of
length $3$ consisting of $x_1,x_2,x_3$. Note that $|L \otimes
\cI_{\xi}|$ is nonempty, since $H^0(X,L)=7$. Moreover, for any $D
\in |L \otimes \cI_{\xi}|$ we have
$$(D \cdot \Gamma)_{x_1} \geq 2, \quad
\text{\rm{and}} \quad D \cap \Gamma \supset \{x_2, x_3\}.$$ Thus
$$4 \leq D \cdot \Gamma = (-3K_X)\cdot(-K_X) = 3.$$
It thus follows that $D$ is reducible and contains $\Gamma$.
Therefore $\xi \in \cB_3^0(X,L)$ and
$b_0(X,L)=3$. Note however that it cannot be $b(X,L)=2$, since
$(X,L)$ does not contain lines. Therefore $b(X,L)=3$. Note also
that the fixed component of $|L \otimes \cI_{\xi}|$ is a singular
plane cubic in the embedding given by $L$.
\end{exmp}

As mentioned in the introduction, the main goal of this work is to shed light
on the global relationship between $n,$ $k$ and $b$ in the case of a $k$-very ample line bundle.
Theorem \ref{dimbound} and Proposition \ref{lowerbound} give lower bounds for $b$ respectively in terms of
$n$ and $k.$ The fact that equality in Proposition \ref{lowerbound} occurs for $n=2,$ see
Remark \ref{p20k}, suggests $b \ge n+k -1$ as a reasonable bound to expect.
The following Theorem proves the suggested bound for $n \ge 4.$ Unfortunately,
 the information given by the same argument for $n=3$ is unsatisfactory.

\begin{thm}
\label{betterbound}
 Let  $L$ be a $k$-very ample line bundle on a
projective $n$-fold $X,$ with $n \ge 3$ and $k \ge 2.$ Assume
there exists $\xi \in \badL{b}$ such that $(H^0(L), \xi)$ is
suitable.
\begin{enumerate}
\item[a)] If $n \ge 4$ then $b(X,L) \ge n + k -1.$ Moreover, if
equality holds, then $|L \otimes \iofo{\xi}|$ has no fixed
component for all $\xi \in \badL{b}$ such that $(H^0(L),\xi)$ is
suitable.
 \item [b)] If $n = 3$ and $b = k+1$ then $k\ge 3.$ Moreover,
 for all $\xi \in \badL{b}$ such that $(H^0(L),\xi)$ is
suitable, $|L \otimes \iofo{\xi}|$ has no fixed component and
every irreducible component of $\bslocus{L \tensor \iofo{\xi}}$ is
a rational normal curve of degree $k.$
\end{enumerate}
\end{thm}
\begin{pf}
As $L$ is very ample, we identify any $0$-scheme $\xi$ with its image in the embedding $\phiL$
and by the linear span of $\xi$ we mean that of $\phiL(\xi).$

From Proposition \ref{lowerbound} it is $b = k+m$ with $m \ge 1.$
Let $\xi \in \badL{b}$ be such that $(H^0(L), \xi)$ is suitable
and let $\Gamma$ be a component of maximal dimension of
$\bslocus{L\tensor \iofo{\xi}}.$ Set $t = \dim{\Gamma}.$ Remark
\ref{pencil} gives $ t = n-1$ or $n-2.$ It is also $\Gamma
\subseteq <\xi> = \Pin{k+m-1}$ because of Lemma \ref{conditions}.
Let $<\Gamma> = \Pin{r}$ so that $t \le r$ and
\begin{equation}
\label{star} r \le k + m -1.
\end{equation}
Assume  that
\begin{equation}
\label{nmassumption2} n \ge m+2,
\end{equation}
or \begin{equation} \label{nmassumption1} n = m+1 \text{  and   }
t = n-1.
\end{equation} In both cases it is:
\begin{equation}
\label{doublestar} t \ge m. \end{equation} As $k \ge 2,$ $\Gamma$
cannot contain lines. Thus if $\Gamma$ is a variety of minimal
degree, i.e. $\deg{\Gamma} =  \text{codim}_{<\Gamma>} \Gamma +1 = r - t +
1,$ then one of the following cases must occur \cite[Theorem 1]{EH}:
\begin{enumerate}
\item[1)] $\Gamma$ is the Veronese surface in $\Pin{5};$
 \item[2)]$ \Gamma$ is a rational
normal curve of degree $r.$
\end{enumerate}
In the former case  \brref{star} implies $m\ge 4$ which contradicts
\brref{doublestar}. In the latter case it is $r\ge k$ as $L$ is
$k$-very ample. This fact combined with \brref{star} and
\brref{doublestar} implies that $r=k,$ and $m=1,$ as $t = 1.$ The
assumption $n\ge 3$ then implies $n=3$ as in part b) of the
statement.

In view of the above argument we can assume that $\deg{\Gamma}  >
r - t + 1.$ Let $\lambda$ be a zero-subscheme of $\Gamma$
consisting of
 $r-t+1$ linearly independent points. Then its linear span
  $\Lambda = <\lambda>$ is a $\Pin{r-t}.$ As $\length{\Lambda \cap \Gamma} = \deg{\Gamma} >\length{\lambda}$ there exists
  a zero-scheme $\lambda'$ of $\Gamma \cap \Lambda$ with $\length {\lambda'} = \length{\lambda} +
  1$ and in particular
  \begin{equation}
  \label{samesystems}
  |L\tensor \iofo{\lambda}| = |L\tensor \iofo{\lambda'}|.
  \end{equation}
   Notice that \brref{star} and \brref{doublestar} imply $r-t+1 \le
   k,$ hence $\length{\lambda'} \le k+1$ contradicting
   the $k$-very ampleness of $L$ in view of \brref{samesystems}.

Consequently neither of the assumptions \brref{nmassumption1} and
\brref{nmassumption2} can hold unless the setting is as in part
$b)$ of the statement, in which case Proposition \ref{b=3} implies
$k \ge 3.$  Therefore if $n
   \ge 4$ and
    $n = m+1 = b-k + 1$ then $t= n-2,$ i.e. $|L \tensor \iofo{\xi}|$ has no fixed component for all $\xi \in \badL{b},$
    such that $(H^0(L),\xi)$ is suitable.
    Moreover, if $n \ge
   4$ it must be $ b-k = m > n -2 $ i. e. $b \ge n + k -1.$  This completes
    the proof of $a).$
\end{pf}
\begin{rem}
\label{surfacesbk+1} If $n = 2,$ and $b = k+1,$ the $k$-very
ampleness of $L$ forces a potential fixed component $\Gamma$ of
$|L \tensor \iofo{\xi}|,$ for all $\xi \in \badL{k+1},$ to be a
rational normal curve of degree $k.$ This can be seen by using a
simple adaptation of the main argument of the proof of Theorem
\ref{betterbound}. There are
no examples known to us of threefolds $(X,L)$ with $b = k+1$ and $k
\ge 3.$
\end{rem}


\section{Reduced bad zero-schemes}

Let $L$ be $k$-very ample, $k\ge 1,$ and assume that
$\badO{t}(X,L)$ is not empty for some $t.$
 On the basis of
\cite{peculiar1} and \cite{DeF}, the naive approach would be to
consider the blow-up $\pi:\tilde{X} \to X$ of $X$ at one point in
the $t$-th bad locus, hoping to obtain, inductively, a new pair,
$(\tX,\pi^*(L)-E),$ polarized with a line bundle which is still
very positive, and has a non empty $\badO{t-1}.$ This is,
unfortunately, not a good strategy. The new polarized pair
contains a linear  $\Pin{n-1},$ which is impossible for a $k$-very
ample line bundle if $k\ge 2.$ Nonetheless, proceeding with a
little care in the same context, it is possible to obtain a new
pair, polarized with a very ample line bundle admitting a non
empty $\badO{2}.$

The following example will shed more light on the situation.
\begin{exmp}
 \label{p2o2} Let $(X,L) = (\Pin{2}, \oofp{2}{2}).$
Notice that $L$ is $2$-very ample.
 Let $\ell$ be
any line in $\Pin{2}$ and let $x_1,x_2,x_3$ be any three collinear
points on $\ell.$ A conic through the $x_i$s must contain $\ell$
and therefore be reducible. Thus the reduced zero-scheme $\xi =
\{x_1, x_2, x_3\}$ is contained in $\badO{3}(X,L).$ This is an
example of a $k$-very ample line bundle with non empty
$\badO{k+1}.$ Notice that $|L \otimes \iofo{\xi}|$ has $\ell$ as
fixed component.

Let now $\pi_1: \tilde{X}_1 \to X$ be the blow-up of $X$ at $x_1$
with exceptional divisor $E_1.$ It is $\mathcal{L}_1 = \pi_1^*(L)
-E_1 = E_1 + 2f$ where $f$ is the proper transform of a line in
$X$ through $x_1.$ Thus ${\rm Pic}(\tilde{X}_1) = \mathbb{Z}[E_1]
\oplus \mathbb{Z}[f].$ Notice that $\cL_1$ is very ample and
embeds $\tilde{X}_1$ in $\Pin{4}$ as a rational normal scroll of
degree $3.$ Let $y_2$ and $y_3$ on $\tilde{X}_1$ be such that
$\pi_1(y_i) = x_i.$ Then $\{ y_2, y_3\} \in \badO{2}(\tilde{X}_1,
\cL_1).$ The situation here is exactly as described in
\cite[Theorem 1.7.9]{BESO} with $\ell = f$ and $R = E_1 + f.$

Let now $\pi_2: \tilde{X}_2 \to X$ be the blow-up of $X$ at $x_1$
and $x_2$ with exceptional divisors $E_1$ and $E_2.$ Consider the
line bundle $\mathcal{L}_2 = \pi_2^*(L) -E_1 - E_2.$ This time
$\cL_2$ is merely spanned and not ample. Let $y_3$ on
$\tilde{X}_2$ be such that $\pi_2(y_3) = x_3.$ Then $ y_3 \in
\cB_1(\tilde{X}_2, \cL_2).$ The divisor $\tilde{\ell} =
\pi_2^*(\ell) - E_1 - E_2$ is a (-1)-curve through $y_3$,
contracted by $\phi_{\cL_2}$ to a smooth point on a smooth quadric
surface. The situation here is exactly as described in
\cite[Theorem 1.1 case d)] {DeF}. Here $\cB_1(\tilde{X}_2, \cL_2)
= \tilde{\ell}.$
\end{exmp}

Pairs $(X,L)$ with $L$ being $k$-very ample, whose $b_0$ index
achieves the lower bound given in Proposition \ref{lowerbound},
turn out to be only $2$-dimensional. As \cite[Theorem 1.7.9]{BESO}
and Proposition \ref{b=3} suggest, one would expect that $L$
$k$-very ample and $b = k+1$ should imply, for all $\xi \in
\badL{b},$
 the existence of a fixed component for $|L \tensor \iofo{\xi}|$ which is a rational
 normal curve of degree $k.$
 The following Proposition gives the desired characterization assuming $b_0 (X,L) =
 k+1,$ and $k \ge 2.$

\begin{prop}
\label{b0=k+1} Let $L$ be a $k$-very ample line bundle on an
$n$-dimensional manifold $X$ with $k \ge 2.$ Assume
$b_0(X,L)=k+1.$ Then $\dim{X} = 2$ and for all $\xi \in
\badO{k+1}$ the linear system $|L \otimes \iofo{\xi}|$ has a
smooth fixed component $\Gamma,$ embedded by $|L|$ as a rational
normal curve of degree $k,$ such that $\xi \subset \Gamma.$
\end{prop}
\begin{pf}
Let $\xi=(x_1,\dots,x_{k+1})\in \cB^0_{k+1}(X,L)$. Let
$\pi:\tilde{X}\to X$ be the blow-up of $X$ at $x_1,\dots,x_{k-1}$
and $y_i=\pi^{-1}(x_{k+1-i})$ for $i=1,2.$ Let $\eta=\{ y_1,
y_2\}.$ Lemma \ref{veryample} implies that $\cL=\pi^*(L)-E_1-\
\dots \ -E_{k-1}$ is very ample. Because $|\cL\otimes
\cI_{\eta}|=|L\otimes\cI_{\xi}|,$ it is $\eta \in
\badO{2}(\tilde{X},\cL)$. Therefore $(\tilde{X},\cL)$ must be as
in \cite[Theorem 1.7.9]{BESO}, i.e. $\dim{X} = 2$ and there exists
a line $\ell$ through $\eta$  with  $|\cL  \otimes \iofo{\eta}| =
\ell + |R|$ where $\oof{X}{R}$ is spanned. Let $\Gamma =
\pi(\ell).$ Because $L$ is $k$-very ample it must be $L\cdot
\Gamma \ge k,$ with equality holding only if $\Gamma \simeq
\Pin{1}.$ Because $\ell$ and $E_i$ are $\cL$-lines, it must be
$\nu_i :=\ell \cdot E_i \in \{ 0,1\}.$ Therefore it is $1 = \cL
\cdot \ell = (\pi^*(\Gamma)- \sum_i \nu_i E_i)\cdot (\pi^*(L) -
\sum_i E_i) = \Gamma \cdot L - \sum_i \nu_i \ge k -(k-1) = 1$
which implies $\Gamma \cdot L = k,$ i.e. $\Gamma$ is a smooth
$\Pin{1},$ embedded by $|L|$ as a rational normal curve of degree
$k.$ Moreover $\ell \cdot E_i = 1$ for all $i,$ i.e. all $x_i$'s
belong to $\Gamma.$
\end{pf}
%
%
Pairs $(X,L)$ where $L$ is $k$-very ample, for which $b_0 = k+2,$
can be further described under the assumption that $\dim{X} = 2.$
The case $k=1$ is analyzed first. The general case is then
obtained from it.
\begin{prop}
\label{b0=3} Let $L$ be a very ample line bundle on a surface $X$.
Assume $b_0(X,L)= 3.$ For all $\xi \in \badO{3},$ if $|L \tensor
\iofo{\xi}|$ does not have  a fixed component then $\phiL$ embeds
$X$ in $\Pin{N},$ in such a way that there exists a linear
$\Pin{2}\subset \Pin{N},$ tangent to $\phiL(X)$ at the $3$ points
$\phiL(\supp{\xi}).$
\end{prop}
\begin{pf}
Due to the assumption on the absence of a fixed component,
$|L\otimes \iofo{\xi}|$ has finite base locus. Hence, by
Proposition \ref{multiplicity}, every $D \in |L \otimes
\iofo{\xi}|$ is of the form $D=A_{b_1}+\ \dots \  +A_{b_r}$, $r
\geq 2$, where all $A_{b_j}$'s belong to a rational pencil $B$.
Moreover, if $\supp{\xi} = \{x_1,x_2,x_3\},$ Proposition \ref{multiplicity}, $iii),$ implies that  $r \geq 2$ and that every $D \in |L\otimes \iofo{\xi}|$ has a
point of multiplicity $\geq 2$ at $x_i$ for $i=1,2,3.$ This gives
the following chain of equalities:
$$\bigcap_{i=1}^{3} |L \otimes \mathfrak{m}_{x_i}^2| =
|L\otimes \iofo{\xi}^2| = |L\otimes \iofo{\xi}|.$$ By Lemma
\ref{conditions} the term on the right is a linear subspace of
codimension $3$ in $|L|$. On the other hand, each of the linear
spaces $|L\otimes \mathfrak{m}_{x_i}^2|$ appearing on the left has
codimension $3$ in $|L|$, since $L$ is very ample, see for example
\cite[Proposition 1.3 and Remark 2.3.3]{la-pa-so2}. This means
that the three linear subspaces $|L\otimes \mathfrak{m}_{x_i}^2|$
coincide, i.\ e.,
$$|L\otimes \mathfrak{m}_{x_1}^2|=|L\otimes \mathfrak{m}_{x_{2}}^2|=|L\otimes
\mathfrak{m}_{x_3}^2|.$$ In other words, looking at $X$ embedded
by $|L|$, every hyperplane tangent to $X$ at $x_1$ is tangent also
at $x_2$ and $x_3$. Equivalently, $X$ embedded by $|L|$ has the
same embedded tangent plane at $x_1, x_2, x_3$.
\end{pf}

\begin{prop}
\label{b0=k+2} Let $k\geq 1$ and $L$ be a $k$-very ample line
bundle on a surface $X.$ Assume $b_0(X,L)=k+2.$ For all $\xi \in
\badO{k+2},$ if $|L \tensor \iofo{\xi}|$ does not have a fixed
component then there exists a linear $\Pin{k+1},$ tangent to
$\phiL(X)$ at the $k+2$ points $\phiL(\supp{\xi}).$
\end{prop}
\begin{pf}
Let $\xi=\{x_1,\dots, x_{k+2}\} \in \badO{k+2}(X,L).$ Let $\eta
\in \red{k-1}$ be any reduced zero-scheme obtained by choosing
$k-1$ of the $k+2$ points of $\xi.$ Let $\tau = \xi \setminus \eta
= \{ x,y,z\}.$  Let $\pi:\tX \to X$ be the blow-up of $X$ at the
$k-1$ points of $\eta$ with exceptional divisors $E_i,$ let $\cL =
\pi^*(L) - \sum_{i=1}^{k-1}E_i,$ and let $\ttau = \{ \tilde{x},
\tilde{y}, \tilde{z}\} = \pi^{-1}(\tau).$ According to Lemma
\ref{veryample}, $\cL$ is very ample on $\tX$ and the proof of
Proposition \ref{b0=3} gives $|\cL \tensor \iofo{\ttau}| = |\cL -
2\tx -2\ty -2\tz|.$ As $|L \tensor \iofo{\eta}| = |\cL|$ we have
$|L \tensor \iofo{\xi}| = |L \tensor \iofo{\eta} \tensor
\iofo{\tau}^2|.$ Because the last equality is true no matter how
$\eta$ was chosen, it follows that $|L \tensor \iofo{\eta}| =
|L\tensor \iofo{\eta}^2|.$ This means that the $\Pin{k+1}$ spanned
by $\phiL(\xi)$ is tangent to $\phiL(X)$ at the $k+2$ points
$\phiL(\supp{\xi}).$
\end{pf}
%
%
\section{Bad Linear Spaces}
As mentioned in the introduction, one may view a bad point as a
{\it bad linear space} of codimension two, as Andrew Sommese
suggested to the first author. In this section we adopt this point
of view. After introducing a natural definition of bad linear
spaces we show that they must necessarily have codimension two and
that they are inherited by hyperplane sections. These two facts
are combined to show that bad linear spaces of very ample linear
systems do not occur at all.

\begin{defn}
Let $X$ be a smooth projective $n$-fold, $n\ge 2,$ and let $L$
 be an ample line bundle on $X$ spanned by $V \subseteq H^0(X,L).$ Let $\Lambda \subset X$
be an $L$-linear subspace of codimension $\geq 2$, i.e. $(\Lambda,
\restrict{L}{\Lambda}) = (\Pin{r}, \oofp{r}{1}),$ for some $r \le
n-2.$ Let $\iofo{\Lambda}$ be the ideal sheaf of $\Lambda.$ We say
that $\Lambda$ is a {\it bad linear space} for $(X,V)$ if for all
$D \in |V \tensor \iofo{\Lambda}|,$ $D$ is reducible or
non-reduced.
\end{defn}
\begin{lem}
\label{codimbadspacesis2}
 Let $X$ be a smooth projective $n$-fold, $n \geq 2,$ and let $L$ be an ample line bundle on $X$,
spanned by a subspace $V \subseteq H^0(X,L)$. If $\Lambda$ is a
bad linear space of $(X,V)$, then $\text{{\rm codim}}_X(\Lambda) =
2$.
\end{lem}
\begin{pf}
Let $\xi$ be a zero-scheme on $X$ consisting of $r+1$ distinct points on $\Lambda,$ not
lying on an
$L$-hyperplane of $\Lambda,$ so that $|V \tensor \iofo{\Lambda}|=|V \tensor \iofo{\xi}|$ with $\xi$
imposing $r+1$ independent conditions on $V.$ Remark \ref{smallBaseLocus} gives
$\dim{(\bslocus{V \tensor \iofo{\xi}})} \le r.$ On the other hand $\Lambda \subset \bslocus{V \tensor \iofo{\xi}},$ hence
$\dim{(\bslocus{V \tensor \iofo{\Lambda}})} = r.$ In particular, as $r\le n-2,$
$|V \tensor \iofo{\Lambda}|$ has
no fixed component.
Thus it follows from Remark \ref{pencil} that $r = \dim{\bslocus{V \tensor \iofo{\Lambda}}}= n-2.$
\end{pf}

The following Proposition shows that bad linear spaces are inherited by hyperplane sections.
To see this, let $\xel$ and $\Lambda$ be as above. Let $x \in \Lambda.$ As $V$ spans $L,$ there exists
a smooth $Y \in |V|$ not passing through $x$ and in particular $\Lambda \not \subset Y.$
Then $\lambda := \Lambda \cap Y$ is an $\restrict{L}{Y}$-hyperplane of $Y.$ Let $W$ be the image
of $V$ under the restriction homomorphism $H^0(X,L) \to H^0(Y,\restrict{L}{Y}).$

\begin{prop}
\label{badhypsec}
 Let notation be as above. If $\Lambda$ is a
bad linear space for $(X,V)$, then $\lambda$ is a bad linear space
for $(Y,W)$.
\end{prop}
\begin{pf}
Let $\rho: V\to W$ be the  homomorphism induced by the restriction
$H^0(X,L) \to H^0(Y,\restrict{L}{Y}).$
Clearly $\rho$ is a surjection and its kernel is $\mathbb{C}
\langle s_0 \rangle$, where $s_0 \in V$ is a non-trivial section vanishing on
$Y$. So, $\dim (W) = \dim (V) - 1.$ Note that $\lambda = \mathbb
P^{n-3},$ by Lemma \ref{codimbadspacesis2} and $\rho(V \otimes \mathcal{I}_{\Lambda})\subseteq W \otimes \mathcal{I}_{\lambda}.$
 Moreover,
$\text{Ker}(\rho) \cap (V \otimes \mathcal{I}_{\Lambda}) = \{ 0
\},$ because any non-trivial element of $\text{Ker}(\rho)$
vanishes exactly on $Y$, hence it cannot vanish on $\Lambda$.
Therefore the homomorphism
\begin{equation}
\label{rhorestr}
\restrict{\rho}{V \otimes \mathcal{I}_{\Lambda}}: V \otimes \mathcal{I}_{\Lambda} \to W
\otimes \mathcal{I}_{\lambda}
\end{equation} is an injection. On the other
hand
\begin{alignat}{1}
\notag \dim (V \otimes \mathcal{I}_{\Lambda}) = &\dim(V) - (\dim (\Lambda) +1)\\ \notag
=&\dim (W) - (\dim (\lambda) +1) \\ \notag = &\dim (W \otimes \mathcal{I}_{\lambda}).
\end{alignat}
 Hence \brref{rhorestr} is an isomorphism, which gives the assertion.
\end{pf}

\begin{thm} Let $X$ be a smooth projective variety of
dimension $n \geq 2$ and let $|V|$ be a very ample linear system
on $X$. Then $(X,V)$ cannot contain bad linear spaces.
\end{thm}
\begin{pf}
Let $S$ be the smooth surface  cut out by $(n-2)$ general elements
of $|V|$ and let $|U|$ be the trace of $|V|$ on $S$. Note that the
corresponding linear system $|U|$ is very ample on $S.$ Now, by
contradiction let $\Lambda$ be a bad linear space for $(X,V).$ By
an inductive application of Proposition \ref{badhypsec} and Lemma
\ref{codimbadspacesis2} we conclude that $p:= \Lambda \cap S$ is a
bad point of $(S,U)$. This contradicts \cite[Corollary 2,
ii)]{peculiar1}.
\end{pf}


\end{document}